\documentclass{gtart}

\def\ifplaintex{\expandafter\ifx\csname documentclass\endcsname\relax}

\ifplaintex 
\else
\fi

\expandafter\ifx\csname beginpicture\endcsname\relax
\expandafter\ifx\csname documentclass\endcsname\relax
\input pictex \else
\input prepictex \input pictex \input postpictex \fi\fi

\def\gt{{\mathsurround=0pt\it $\cal G\mskip-2mu$eometry \&\ 
$\cal T\!\!$opology}}        

\def\gtp{{\mathsurround=0pt\it $\cal G\mskip-2mu$eometry \&\ 
$\cal T\!\!$opology $\cal P\!$ublications}}  

\def\lognumber#1{\def\thelognumber{#1}}
\def\volumenumber#1{\def\thevolumenumber{#1}}
\def\papernumber#1{\def\thepapernumber{#1}}
\def\volumeyear#1{\def\thevolumeyear{#1}}

\def\pagenumbers#1#2{\def\startpage{#1}\def\finishpage{#2}}
\def\published#1{\def\publishdate{#1}}
\def\proposed#1{\def\theproposer{#1}}
\def\seconded#1{\def\theseconders{#1}}
\def\received#1{\def\receiveddate{#1}}
\def\revised#1{\def\reviseddate{#1}}
\def\accepted#1{\def\accepteddate{#1}}

\long\def\asciiabstract#1{\long\def\theasciiabstract{#1}}
\def\asciikeywords#1{\def\theasciikeywords{#1}}

\let\\\par\let\thelognumber\relax
\let\thevolumenumber\relax\let\thepapernumber\relax
\let\thevolumeyear\relax\let\thesamplenumber\relax\let\startpage\relax
\let\finishpage\relax\let\publishdate\relax\let\receiveddate\relax
\let\reviseddate\relax\let\accepteddate\relax\let\theasciititle\relax
\let\theasciiauthors\relax
\let\theasciiabstract\relax\let\theasciikeywords\relax
\let\theasciiemail\relax\let\theshortauthors\relax\let\theshorttitle\relax

\long\def\maketitlep{   

\count0=\startpage

\gt\hfill      
\beginpicture
\setcoordinatesystem units <0.33truein, 0.33truein> point at 2.2 0.9
\setplotsymbol ({$\cal G$})
\plotsymbolspacing=9truept
\circulararc 315 degrees from 0 1 center at 0 0
\setplotsymbol ({$\cal T$})
\circulararc 315 degrees from 1 -1 center at 1 0
\endpicture
\break
{\small\ifx\thesamplenumber\relax 
Volume \else Sample
\fi\thevolumenumber\ (\thevolumeyear)
\startpage--\finishpage\nl
Published: \publishdate}
\vglue 0.5truein plus 0.4fil minus 0.1truein

{\parskip=0pt\leftskip 0pt plus 1fil\def\\{\par\smallskip}{\ifplaintex\large
\else\Large\fi\bf\thetitle}\par\medskip}   

\vglue 0pt plus 0.1fil 

{\parskip=0pt\leftskip 0pt plus 1fil\def\\{\par}{\sc\theauthors}
\par\medskip}

\vglue 0pt plus 0.1fil 

{\small\parskip=0pt\let\newline\\
{\leftskip 0pt plus 1fil\def\\{\par}{\sl\theaddress}\par}
\expandafter\ifx\theemail\relax    
\relax\else\vglue 5pt plus 0.02fil minus 2pt\def\\{\stdspace{\rm 
and}\stdspace} 
\cl{Email:\stdspace\tt\theemail}\fi
\ifx\theurl\relax                  
\relax\else\vglue 5pt plus 0.02fil minus 2pt\def\\{\stdspace{\rm 
and}\stdspace}
\cl{URL:\stdspace\tt\theurl}\fi\par}

\vglue 7pt plus 0.3fil minus 3pt

{\bf Abstract}
\vglue 5pt plus 0.1fil minus 2pt

\theabstract

\vglue 7pt plus 0.3fil minus 3pt

{\bf AMS Classification numbers}\quad Primary:\quad \theprimaryclass

Secondary:\quad \thesecondaryclass

\vglue 5pt plus 0.3fil minus 2pt

{\bf Keywords}\quad \thekeywords

\vglue 10pt plus 0.5fil minus 5pt

{\small  Proposed: \theproposer\hfill Received: \receiveddate\nl
Seconded: \theseconders\hfill 
\ifx\reviseddate\relax                         
Accepted: \accepteddate                        
\else
Revised: \reviseddate                          
\fi}
\eject
}       

\let\maketitlepage\maketitlep
\let\maketitle\maketitlepage

\font\phead=cmsl9 scaled 950
\font=cmsl9 scaled 1050
\font\pnum=cmbx10 scaled 913
\font\lnum=cmbx10 
\font\pfoot=cmsl9 scaled 950
\font=cmsl9 scaled 1050
\ifplaintex
\headline{\vbox to 0pt{\vskip -4.5mm\line{\small\phead\ifnum
\count0=\startpage ISSN 1364-0380 (on line)
1465-3060 (printed) \hfill {\pnum\folio}\else\ifodd\count0\def\\{ }%
\ifx\theshorttitle\relax\thetitle\else\theshorttitle\fi\hfill{\pnum\folio}
\else\def\\{ and }{\pnum\folio}\hfill\ifx\theshortauthors\relax\theauthors
\else\theshortauthors\fi\fi\fi}\vss}}
\footline{\vbox to 0pt{\vglue 0mm\line{\small\pfoot\ifnum\count0=\startpage
\copyright\ \gtp\hfill\else
\gt, Volume \thevolumenumber\ (\thevolumeyear)\hfill\fi}\vss
}}
\else
\makeatletter
\def\@oddhead{{\small\lhead\ifnum\count0=\startpage ISSN 1364-0380 (on line)
1465-3060 (printed) \hfill {\lnum\number\count0}\else\ifodd\count0
\def\\{ }\ifx\theshorttitle\relax \thetitle \else\theshorttitle\fi\hfill
{\lnum\number\count0}\else\def\\{ and }{\lnum\number\count0}
\hfill\ifx\theshortauthors\relax 
\theauthors\else\theshortauthors\fi\fi\fi}}\def\@evenhead{\@oddhead}
\def\@oddfoot{\small\lfoot\ifnum\count0=\startpage\copyright\ \gtp\hfill\else
\gt, Volume \thevolumenumber\ (\thevolumeyear)\hfill\fi}
\def\@evenfoot{\@oddfoot}
\makeatother
\fi

\newwrite\gtoutfile
\long\gdef\makeheadfile{  
{\def\\{, }\def\s{ }
\immediate\openout\gtoutfile head.xxx
\immediate\write\gtoutfile{To: math@arxiv.org}
\immediate\write\gtoutfile{Subject: put or rep NNNNN:pppp}
\immediate\write\gtoutfile{--text follows this line--}
\immediate\write\gtoutfile{Proxy-for: \ifx\theasciiauthors\relax
\theauthors\else\theasciiauthors\fi\s<\ifx\theasciiemail\relax\theemail\else\theasciiemail\fi>}
\immediate\write\gtoutfile{\noexpand\\}
\immediate\write\gtoutfile{Authors: \ifx\theasciiauthors\relax
\theauthors\else\theasciiauthors\fi}
{\def\\{ }\immediate\write\gtoutfile{Title: \ifx\theasciititle\relax
\thetitle\else\theasciititle\fi}}
\immediate\write\gtoutfile{Subj-class: GT or SG or MG etc}
\immediate\write\gtoutfile{MSC-class: \theprimaryclass\ifx\thesecondaryclass\relax\else, \thesecondaryclass\fi}
\immediate\write\gtoutfile{Journal-ref: Geom. Topol. \thevolumenumber
(\thevolumeyear) \startpage-\finishpage}
\immediate\write\gtoutfile{Comments: Published by Geometry and Topology at}
\immediate\write\gtoutfile{\s\s http://www.maths.warwick.ac.uk/gt/GTVol\thevolumenumber/paper\thepapernumber.abs.html}
\immediate\write\gtoutfile{\noexpand\\}
\immediate\write\gtoutfile{}
\ifx\theasciiabstract\relax
\immediate\write\gtoutfile{\theabstract}\else
\immediate\write\gtoutfile{\theasciiabstract}\fi
\immediate\write\gtoutfile{}
\immediate\write\gtoutfile{\noexpand\\}
\immediate\write\gtoutfile{}
\immediate\closeout\gtoutfile}}  

\def\maketitlepage{\maketitlep\makeheadfile}
\let\maketitle\maketitlepage

\lognumber{162}

\volumenumber{6}
\papernumber{27} 
\volumeyear{2002}
\pagenumbers{905}{916} 
\received{17 January 2001}
\revised{7 November 2002}
\accepted{31 December 2002}
\published{31 December 2002}
\proposed{David Gabai}
\seconded{Walter Neumann, John Morgan}

\usepackage{amsmath,amssymb,epsfig,psfrag,subfigure}

\newtheorem{theorem}{Theorem}[section]
\newtheorem{lemma}[theorem]{Lemma}
\newtheorem{corollary}[theorem]{Corollary}

\theoremstyle{definition}

\newcommand{\RR}{\mathbb{R}}

\def\vol{\mbox{\rm{Vol}}}

\def\a{\alpha}
\def\b{\beta}
\def\e{\epsilon}
\def\g{\gamma}

\def\th{\theta}
\def\l{\lambda}
\def\d{\delta}
\def\k{\kappa}

\def\i{\iota}

\def\w{\omega}

\def\Vol{\mbox{\rm{Vol}}}
\def\min{\mbox{\rm{min}}}

\begin{document}
\title{Volume change under drilling}
\authors{Ian Agol}
\address{MSCS, SEO 322, m/c 249, University of Illinois at 
Chicago\\851 S Morgan St, Chicago, IL 60607-7045, USA}
\email{agol@math.uic.edu} 
\url{http://www.math.uic.edu/\char'176agol}
\begin{abstract}   
Given a hyperbolic 3--manifold $M$ containing an embedded closed
geodesic, we estimate the volume of a complete hyperbolic metric
on the complement of the geodesic in  terms of the geometry of
$M$. As a corollary, we show that the smallest volume orientable
hyperbolic 3--manifold has volume $>.32$.
\end{abstract}

\asciiabstract{  
Given a hyperbolic 3-manifold M containing an embedded closed
geodesic, we estimate the volume of a complete hyperbolic metric
on the complement of the geodesic in  terms of the geometry of
M. As a corollary, we show that the smallest volume orientable
hyperbolic 3-manifold has volume >.32 .}

\asciikeywords{Hyperbolic structure, 3-manifold, volume, geodesic}

\maketitlepage


\section{Introduction}
 Kerckhoff showed that given a hyperbolic
3--manifold $M$ and an embedded geodesic $\gamma \subset M$, the
manifold $M-\g\equiv M_\g$ has a metric of  negative sectional
curvature \cite{Ko}, using a method of Gromov and Thurston
\cite{BH}. Therefore by Thurston, $M_\gamma$ admits a unique
complete hyperbolic metric, ie, a metric of constant
sectional curvature $-1$  \cite{Mo}. In theorem \ref{bound} of
this paper, we show that if $\gamma$ has an embedded tubular
neighborhood of radius $R$, then $\vol(M_\gamma)\leq (\coth
R)^{\frac52}(\coth 2R)^{\frac12}\ \vol(M)$, where Vol denotes the
volume of the unique complete hyperbolic structure. Specifically,
we use Kerckhoff's method to fill in a metric on $M_\gamma$ in
order to apply a theorem of Boland, Connell, and Souto \cite{BCS},
which is based on an important technique of Besson, Courtois, and
Gallot \cite{BCG}, to get an upper bound on the hyperbolic volume
of $M_\gamma$ in terms of the geometry of $M$. As a corollary, we
show that the smallest volume orientable hyperbolic 3--manifold has
volume $>.32$. It is conjectured that the Weeks manifold, which
has volume $=.9427...$, is the smallest volume hyperbolic
3--manifold. This is known to be the smallest volume arithmetic
hyperbolic 3--manifold, by a result of Chinburg, Friedmann, Jones,
and Reid \cite{CFJR}. The smallest volume orientable cusped
manifolds are known to be the figure eight knot complement and its
sibling, which have volume $= 2.0298...$, by a result of Cao and
Meyerhoff \cite{CM}, which we use combined with the main theorem
to obtain our lower bound.

\medskip
{\bf Acknowledgements}\qua We thank Bill Thurston for helpful
conversations, in particular for explaining lemma \ref{curvature}.
This research was mostly carried out at the University of
Melbourne, and we thank Iain Aitchison, Hyam Rubinstein and the
maths department for their hospitality. We thank Oliver Goodman
for help using his program {\it tube} \cite{Go1}. We thank the
referees for many useful comments, and for finding a mistake in an
earlier version. The author was partially supported by ARC grant
420998.

\section{Estimates}

\begin{theorem}\label{bound}
Let $M$ be an orientable hyperbolic 3--manifold with metric $\nu$.
Let $\gamma \subset M$ be a geodesic in $M$ of length $l$ with an
embedded open tubular neighborhood $C$ of radius $R$, and with
complete hyperbolic metric $\tau$ on $M_\gamma$. 
\begin{gather*}
\text{Then}\qua\vol(M_\gamma,\tau) \leq (\coth R \coth 2R)^{\frac32}
(\Vol(M,\nu)+\pi l \sinh^2R(\frac{\coth R}{\coth 2R}-1))\phantom{TT}\\
\leq (\coth R)^{\frac52}(\coth 2R)^{\frac12}Vol(M,\nu).
\end{gather*}
\end{theorem}

{\bf Remark}\qua This theorem can be easily generalized to a geodesic
link in $M$ with an embedded tubular neighborhood, but we restrict
to the case of one component to simplify the argument, since this
is the context of our main application. Also, if $M$ is
non-orientable, then we may apply this theorem to the 2--fold
orientable covering to recover the same estimate.

The intuition behind the estimate is that we construct a
Riemannian metric $\rho$ on $M_\g$ which coincides with the metric
$\nu$ on $M\backslash C$, such that $Ric_{\rho} \geq -2k \rho$,
where $k=\coth R \coth 2R$. Then we would like to apply the
following Volume Theorem of Boland, Connell and Souto  to compare
the two metrics on $M_\g$:

\begin{theorem}\label{bcs}{\rm\cite{BCS}}\qua 
Let $(M,g)$ and $(M_0,g_0)$ be two oriented complete finite volume
Riemannian manifolds of the same dimension $\geq 3$, and suppose
that $Ric_g\geq -(n-1) g$, and $-a\leq K_{g_0}\leq -1$. Then for
all proper continuous maps $f:M\to M_0$,
$$\Vol(M,g)\geq |deg(f)| \Vol(M_0,g_0)$$
and equality holds if and only if $f$ is proper homotopic to a
Riemannian covering.
\end{theorem}

In our situation, we want to choose $(M_0,g_0)=(M_\g,\tau)$, where
$K_\tau=-1$, since the metric is hyperbolic, $(M,g)=(M_\g,k\rho)$,
where we scale $\rho$ by $k$ so that it satisfies $Ric_{k\rho}\geq
-2k\rho$ almost everywhere, and $f=Id_{M_\g}$ is degree one and
proper. Then we would like to apply theorem \ref{bcs} to conclude
that $\vol(M_\g,k \rho)\geq \vol(M_\g,\tau)$, which may be
computed to obtain the estimate of theorem \ref{bound}. Some
technicalities are involved in that the metric $\rho$ we construct
is only $C^1$, so that we must approximate $\rho$ by smooth
metrics in order to apply theorem \ref{bcs}.

\proof[Proof of theorem \ref{bound}]
The tubular neighborhood $C$ of radius $R$ of $\g$ of length $l$
has a metric in cylindrical coordinates $(r,\th,\l)$ given by
$ds^2=dr^2+\sinh^2 r d\th^2 +\cosh^2 r d\l^2$, on the cylinder
$C=[0,R]\times[0,2\pi]\times[0,l]$ with $(r,0,\l)\sim
(r,2\pi,\l)$, $(r,\th,0)\sim (r,\th+\phi,l)$ and
$(0,\th,\l)\sim(0,0,\l)$, for some constant $\phi$ giving the
rotational part of the holonomy of $\g$. $M_\g$ may be obtained
from $M$ by removing $C$ from $M$ and replacing it with the region
$C'=(-\infty,R]\times[0,2\pi]\times[0,l]$, with the same first two
identifications as above. We will give $C'$  a metric of the form
$ds^2=dr^2+f(r)^2d\th^2+g(r)^2d\l^2$. These types of metrics were
considered by Gromov and Thurston in proving their $2\pi$--theorem
for obtaining negatively curved metrics on Dehn filled 3--manifolds
\cite{BH}.  We need to compute sectional curvatures for this
metric. This has been done by Bleiler and Hodgson \cite{BH}, but
we will give a sketch of a derivation as explained to us by
Thurston.

\begin{lemma}\label{curvature}
The metric $dr^2+f(r)^2d\th^2+g(r)^2d\l^2$ on $\RR^3$
has sectional curvatures
$$K_{\th\l}=-\frac{f'g'}{fg},\ K_{r\th}=-\frac{f''}{f},\
K_{r\l}=-\frac{g''}{g}$$
\end{lemma}
\begin{proof}
First we will compute $K_{\th\l}$. Fix a reference
point $(r_0,\th_0,\l_0)$ at which to compute the
curvature. Consider the plane
$H={(r_0,\th,\l)}$. $H$ is isometric
to the euclidean plane, so the Gauss
curvature of $H$ is 0. The reflection through the plane
$\th=\th_0$
 preserves $H$ and $(r_0,\th_0,\l_0)$,
so the lines of principal curvature of $H$ at
$(r_0,\th_0,\l_0)$ are invariant under this reflection.
Thus, they must be
parallel to the vectors $\frac{\partial}{\partial\th}$
and $\frac{\partial}{\partial\l}$. To compute the
principal curvatures, we can compute the geodesic
curvatures of $H$--geodesics in the $\th$ and $\l$
directions. Consider the curve
$c=(r_0,\th,\l_0)$. The length of $c$ on the interval
$[\th_0-\e,\th_0+\e]$ is $2\e f(r_0)$. If we do a
variation of this curve normal to $H$, the length
is $2\e f(r)$. Since the reflection in the plane $\l=\l_0$
preserves $c$, the osculating plane of $c$ must be the plane
$\l=\l_0$, so that the normal variation lies
in this osculating plane. The geodesic curvature is then the
logarithmic derivative of the length of the variation (compare
to the euclidean case),
$k_g=\frac{\partial}{\partial r} \log (2 \e f(r)) =
\frac{f'}{f}$, which gives the principal curvature
of $H$ in the direction of $\th$. The other principal
curvature is $\frac{g'}{g}$. Thus, the Gauss curvature
of $H$ is $0=K_{\th\l}+ \frac{f'g'}{fg}$, so
$K_{\th\l}=-\frac{f'g'}{fg}$.

To compute $K_{r\th}$, notice that the plane $\l=\l_0$
is geodesic. So the sectional curvature is the
Gauss curvature of the plane $\l=\l_0$.
Consider the annulus $S_r=[r_0,r]\times[0,2\pi]\times\{\l_0\}$.
By the Gauss--Bonnet theorem,
$$\int_{S_r} K_{r\th}\ dA + \int_{\partial S_r} k_g\ ds = 2\pi \chi(S_r)=0.$$
We have
$$\frac{\partial}{\partial r} \int_{S_r} K_{r\th}\ dA=
\frac{\partial}{\partial r} \int_{r_0}^r \int_0^{2\pi} K_{r\th} f(r)\ d\th dr
=
2\pi K_{r\th} f(r)$$

$$= -\frac{\partial}{\partial r}\int_{\partial S_r} k_g\ ds
=-\frac{\partial}{\partial r} \int_0^{2\pi} \frac{f'(r)}{f(r)} f(r)d\th = -2\pi f''(r).$$
Thus $K_{r\th}=-\frac{f''}{f}$. Similarly, $K_{r\l}=-\frac{g''}{g}$.
\end{proof}

\begin{figure}[htb]\small
    \cl{\psfrag{f}{$f$}
    \psfrag{g}{$g$}
    \psfrag{R}{$R$}
    \psfrag{ch}{$\cosh r$}
    \psfrag{sh}{$\sinh r$}
    \epsfxsize=2.5in
    \epsfbox{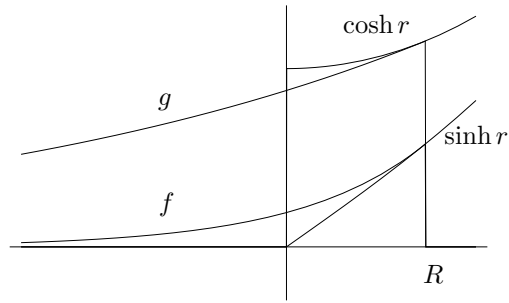}}
    \caption{\label{neg} Constructing a negatively curved metric}
    \end{figure}

Kerckhoff's proof that $M_\g$ has a metric of negative sectional
curvature follows from lemma \ref{curvature} by noting that one
may choose $C^{\infty}$ functions $f(r)$ and $g(r)$ on
$(-\infty,R]$ such that for $r$ near $R$, $f(r)=\sinh r$ and
$g(r)=\cosh r$, and such that all of $f(r), f'(r), f''(r), g(r),
g'(r), g''(r)$ are $>0$ (see figure \ref{neg}). Then it follows
that $M_\g$ is atoroidal, anannular, and irreducible, so $M_\g$
has a complete hyperbolic metric $\tau$  of finite volume by
Thurston's hyperbolization theorem of Haken manifolds \cite{Mo}.

We want to construct metrics on $M_\g$ for which we may apply
theorem \ref{bcs}. Let $f(r)=(\sinh R) e^{ (\coth R )(r-R)}$, and
$g(r)=(\cosh R) e^{(\tanh R) (r-R)}$. This gives a $C^\infty$
metric on $C'$, which induces a $C^1$ metric $\rho$ on $M_\g$,
where $\rho$ agrees with $\nu$ on $M_\g\backslash C'=M\backslash
C$.

On $C'$, choose the orthonormal basis for the tangent space at a
fixed point $(e_1,e_2,e_3)=(\frac{\partial}{\partial
r},\frac{1}{f}\frac{\partial}{\partial\th},\frac{1}{g}\frac{\partial}{\partial\l})$.
The reflections $\theta\to -\theta$ and $\l \to -\l$ preserve the
metric locally, so must preserve the maximal and minimal sectional
curvatures. The only tangent planes preserved by both these
reflections are spanned by $\{e_i,e_j\}$, $i\neq j$. Thus, we
conclude that the sectional curvatures $K_{r\th}=K_{e_1e_2}$,
$K_{\th\l}=K_{e_2e_3}$, and $K_{r\l}=K_{e_1e_3}$ realize the
minimum and maximum sectional curvatures. It is well known that in
such coordinates, the Ricci curvature is diagonalized. To see
this, one computes the curvature tensor $R(e_i,e_j,e_k,e_l)$ (see
\cite{GHL}, ch. III A, for notational conventions). The sectional
curvature of the plane spanned by $\{e_i,e_j\}$ is $R(e_i,e_j
,e_i, e_j)$. If this sectional curvature is maximal, then we see
that
$$\frac{d}{d\varphi} R(e_i,(e_j\cos\varphi + e_k\sin\varphi) , e_i, (e_j\cos\varphi +
e_k\sin\varphi))|_{\varphi=0} = 2 R(e_i, e_k , e_i, e_j) = 0.$$
Similarly, one may see that
$R(e_j,e_k,e_j,e_i)=R(e_k,e_i,e_k,e_j)=0$, by exchanging the roles
of $i$ and $j$, and using the same observation for the minimal
sectional curvature. So all the off-diagonal terms of $R$ vanish.
Since $Ric_\rho(e_i,e_j)=\sum_{k=1}^3 R(e_k,e_i, e_k, e_j)=0$, if
$i\neq j$ we see that the off-diagonal terms of the Ricci tensor
vanish as well. We compute $Ric_\rho$ in the basis
$(e_1,e_2,e_3)$:

 $$ Ric_\rho = \begin{bmatrix} K_{r\th}+K_{r\l} & 0 & 0
 \\ 0 & K_{r\th}+K_{\th\l} & 0
 \\ 0 & 0 & K_{r\l}+K_{\th\l} \end{bmatrix}$$
$$ = \begin{bmatrix} -\coth^2 R -\tanh^2 R & 0 & 0
 \\ 0 & -\coth^2 R - 1 & 0
 \\ 0 & 0 & -\tanh^2 R - 1 \end{bmatrix}$$
 $$\geq -2\coth R\coth 2R
 \begin{bmatrix}1&0&0\\0&1&0\\0&0&1\end{bmatrix}=(-2\coth R\coth 2R) \rho $$
by lemma \ref{curvature}. We have
$$\Vol(M_\g,\rho)=\Vol(M,\nu)- \Vol(C,\nu)+\Vol(C',\rho)$$
$$=\Vol(M,\nu)-\pi l\sinh^2 R + 2\pi l \frac{\sinh R \cosh R}{\coth
R + \tanh R}$$ $$ = \Vol(M,\nu) + \pi l \sinh^2 R (\frac{\coth
R}{\coth 2R} -1) \leq \Vol(M,\nu) \frac{\coth R}{\coth 2R}$$ where
the last inequality uses the fact that $\pi l \sinh^2 R =
\Vol(C,\nu) \leq \Vol(M,\nu)$.

\begin{lemma}\label{sequence}
There is a family of negatively curved metrics $\rho_\e$ on
$M_\g$ such that $\Vol(M_\g,\rho_\e)\to \Vol(M_\g,\rho)$ and
$Ric_{\rho_\e}\geq -2 k_\e \rho_\e$, where $k_\e\to \coth R\coth
2R$ as $\e\to 0^+$.
\end{lemma}

Using this lemma, we can finish the proof of theorem \ref{bound}.
Consider the family of metrics $\rho_\e$ on $M_\g$. The metrics
$k_\e \rho_\e$ on  $M_\g$ satisfy $Ric_{k_\e
\rho_\e}=Ric_{\rho_\e} \geq -2 k_\e \rho_\e$. Given this Ricci
curvature control, we can apply theorem \ref{bcs} with
$(M,g)=(M_\g,\tau)$, $(M_0,g_0)=(M_\g,k_e\rho_\e)$, and
$f=Id_{M_\g}$ to conclude that
$$k_\e^{3/2}\vol(M_\g,\rho_\e) = \vol(M_\g,k_\e \rho_\e)\geq
\vol(M_\g,\tau).$$ In the limit, we have
$$ \vol(M_\g,\tau)\leq (\coth R\coth2R)^{3/2}\vol(M_\g,\rho)$$
$$=(\coth R\coth 2R)^{\frac32}
(\Vol(M,\nu)+\pi l \sinh^2R(\frac{\coth R}{\coth 2R}-1))$$
$$\leq (\coth
R)^{5/2}(\coth 2R)^{1/2}\Vol(M,\nu).\eqno{\qed}$$

To complete the proof of theorem \ref{bound}, we need to prove
lemma \ref{sequence}.

\begin{proof}[Proof of lemma \ref{sequence}]
We will construct the metric $\rho_\e$ by smoothing $\rho$ on
the region $C'$, by smoothing the functions $f(r)$ and $g(r)$
using the construction given in the following lemma.

\begin{lemma}\label{smooth}
Suppose
\begin{equation}
a(r)=
    \begin{cases}
        b(r)& r<R ,\\
        c(r) & r\geq R,
    \end{cases}
\end{equation}
where $b(r)$ and $c(r)$ are $C^\infty$ on $(-\infty,\infty)$, and
$b(R)=c(R)$, $b'(R)=c'(R)$.
Then we may find $C^{\infty}$ functions $a_\e$  on $(-\infty,\infty)$
for $\e>0$ such that
\begin{enumerate}
\item
there is a $\d(\e)>0$ such that $\underset{\e\to 0^+}{\lim}
\d(\e)=0$ and $a_\e(r)=b(r)$ for $r\leq R-\d(e)$, $a_\e(r)=c(r)$
for $r\geq R$, and
\item
$$\min\{b''(R),c''(R)\}=\underset{\e\to 0^+}{\lim} \underset{R-\d(\e)\leq r \leq R}{\inf} a_\e''(r)$$
$$\leq \underset{\e\to 0^+}{\lim} \underset{R-\d(\e)\leq r \leq R}{\sup} a_\e''(r)
= \max\{b''(R),c''(R)\}.$$
\end{enumerate}
\end{lemma}
\begin{proof}[Proof of lemma \ref{smooth}]
The idea is to smoothly interpolate between the functions $b(r)$
and $c(r)$. If we were to just use a partition of unity to
interpolate between them, we would not control the first and
second derivatives of the interpolation. So instead, we use a
partition of unity to interpolate between $b''(r)$ and $c''(r)$
to get a function $\eta_\e'(r)$. Integrating twice, we must add
bump functions in order to get the correct properties.

Define \begin{equation}
\a(r)=
    \begin{cases}
        e^{-1/r^2} e^{-1/(1-r)^2} & 0<r<1 ,\\
        0 & r\leq 0,\ r\geq 1.
    \end{cases}
\end{equation}
Then $\a(r)$ is a $C^\infty$ bump function. Let
 $\b(r)=\int_{0}^r \a(t)dt/\int_0^1 \a(t)dt$.
 Then $\b'(r)\geq 0$, and $\b'$ is supported
 on $[0,1]$.
 Let $\phi_\e(r)=\b((r-R)/\e+1)$, if $\e >0$, and $\phi_0(r)=0$. Let
 $\eta_\e'(r)=b''(r)(1-\phi_\e(r))+c''(r)\phi_\e(r)$ (figure \ref{2nd}(a)). Then
 $\eta_\e'(r)=b''(r)$ for $r\leq R-\e$, and $\eta_\e'(r)=c''(r)$ for $r\geq R$,
 and thus $\eta_\e'(r)$ interpolates between $b''(r)$ and $c''(r)$.
 Let $\eta_\e(r)=b'(R-\e)+\int_{R-\e}^r \eta_\e'(t) dt$.
 $\eta_\e(r)=b'(r)$ for $r\leq R-\e$, but it will differ
 from $c'(r)$ for $r\geq R$ by a constant (figure \ref{deriv}(a)), so we must
 adjust $\eta_\e(r)$ to obtain a function which interpolates between
 $b'$ and $c'$.
 Let $\i(\e)=|c'(R)-\eta_\e(R)|^{1/2}$, then $\lim_{\e\to 0^+} \i(\e)=0$.
 Let $\k_\e'(r) = \eta_\e(r)+ (c'(R)-\eta_\e(R))\phi_{\i(\e)}(r)$ (figure \ref{deriv}(a)).
 Now $\k_\e'(r)= c'(r)$ for $r\geq R$, and $\k_\e'(r)\to a'(r)$
 uniformly as $\e\to 0^+$. We have
  $|c'(R)-\eta_\e(R)|\phi_{\i(\e)}'(r) = \i(\e) \b'((r-R)/\i(\e)+1),$
 which goes to 0 uniformly as $\e\to 0^+$, since $\b'$ is uniformly bounded.
 Let $\k_\e(r)=b(R-\i(\e))+\int_{R-\i(\e)}^r \k_\e'(t) dt$, and define
 $\w(\e)=|c(R)-\k_\e(R)|^{1/3}$. Let
 $a_\e(r)=\k_\e(r)+(c(R)-\k_\e(R))\phi_{\w(\e)}(r)$ (figure \ref{deriv}(b)).
We choose  $\w(\e)$
 so that the bump function has small enough second derivative, that is
$|c(R)-\k_\e(R)|\phi_{\w(\e)}''(r)=\w(\e) \b''((r-R)/\w(\e)+1) \to
0$ as $\e\to 0^+$. Let $\d(\e)=\max\{\e,\i(\e),\w(\e)\}$.  It is
clear from the construction that $a_\e(r)=b(r)$, for $r\leq
R-\d(\e)$, and $a_\e(r)=c(r)$, for $r\geq R$. We also have
$\min\{b''(R),c''(R)\}= \underset{\e\to 0^+}{\lim}
\underset{R-\d(\e)<r<R}{\inf}\eta_\e'(r)\leq \underset{\e\to
0^+}{\lim} \underset{R-\d(\e)<r<R}{\sup}\eta_\e'(r)=
\max\{b''(R),c''(R)\}$. Now, since $|a_\e''(r)-\eta_\e'(r)|\leq
|\i(\e)\b'((r-R)/\i(\e)+1)|+|\w(\e)\b''((r-R)/\w(\e)+1)| \to 0$ as
$\e\to 0^+$ (figure \ref{2nd}(b)), we see that $a_\e''(r)$ satisfies
the second condition of the lemma.
\end{proof}
\begin{figure}[htb]
    \begin{center}
    \begin{huge}
    \psfrag{a}{$a_\e''$}
    \psfrag{d}{$a_{\e/2}''$}
    \psfrag{e}{$a_{\e/4}''$}
    \psfrag{f}{$a_{\e/8}''$}
    \psfrag{x}{$\eta_\e'$}
    \psfrag{y}{$\eta_{\e/2}'$}
    \psfrag{z}{$\eta_{\e/4}'$}
    \psfrag{w}{$\eta_{\e/8}'$}
    \psfrag{b}{$b''$}
    \psfrag{c}{$c''$}
    \psfrag{R}{\raisebox{-2ex}{$R$}}
    \psfrag{W}{\raisebox{-2ex}{$R-\e$}}
    \psfrag{V}{\raisebox{-2ex}{$R-\d(\e)$}}
    \psfrag{U}{}
    \psfrag{T}{}
    \psfrag{S}{}
    \subfigure{\epsfig{figure=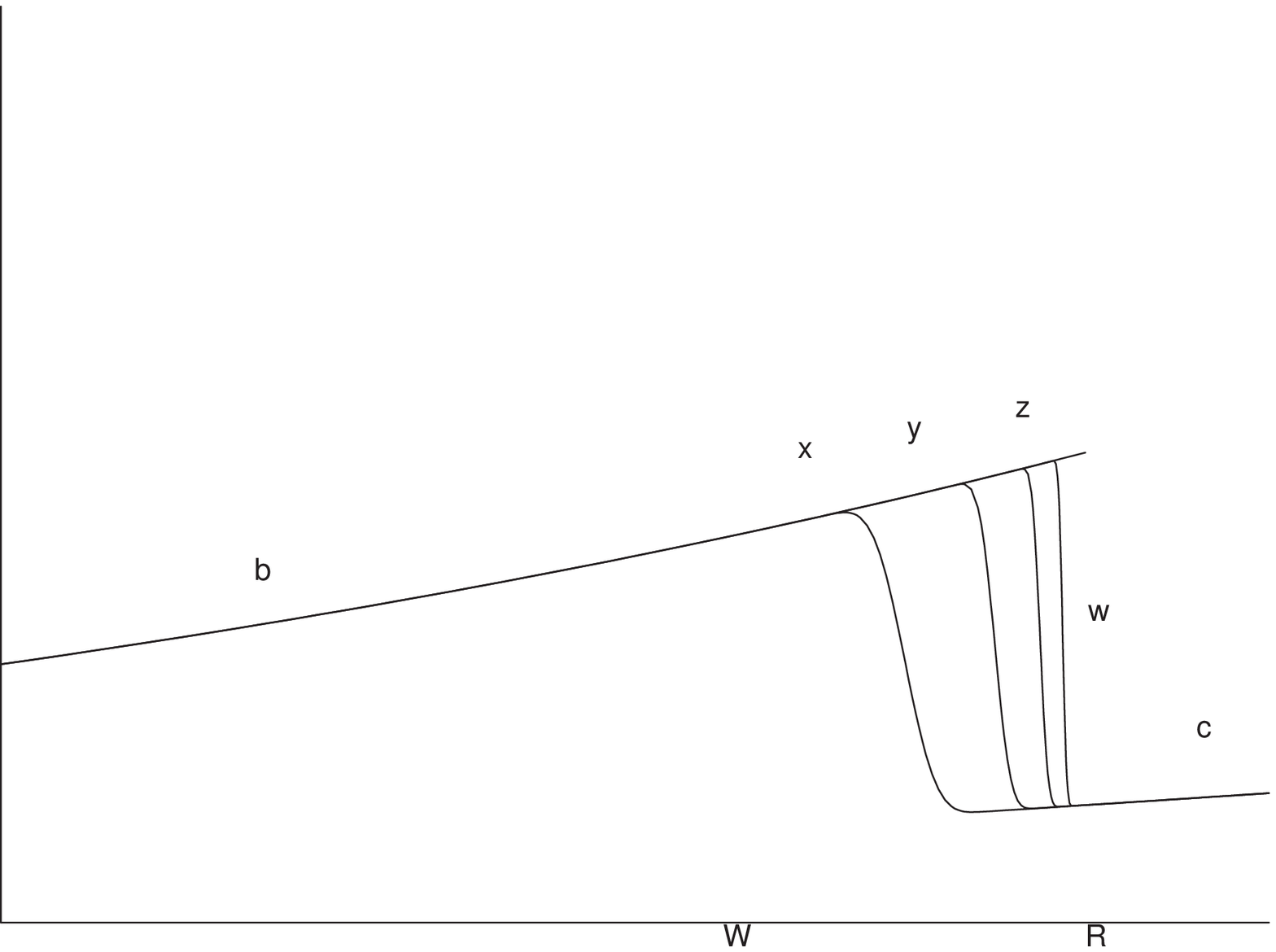,angle=0,width=.45\textwidth}}\quad
    \subfigure{\epsfig{figure=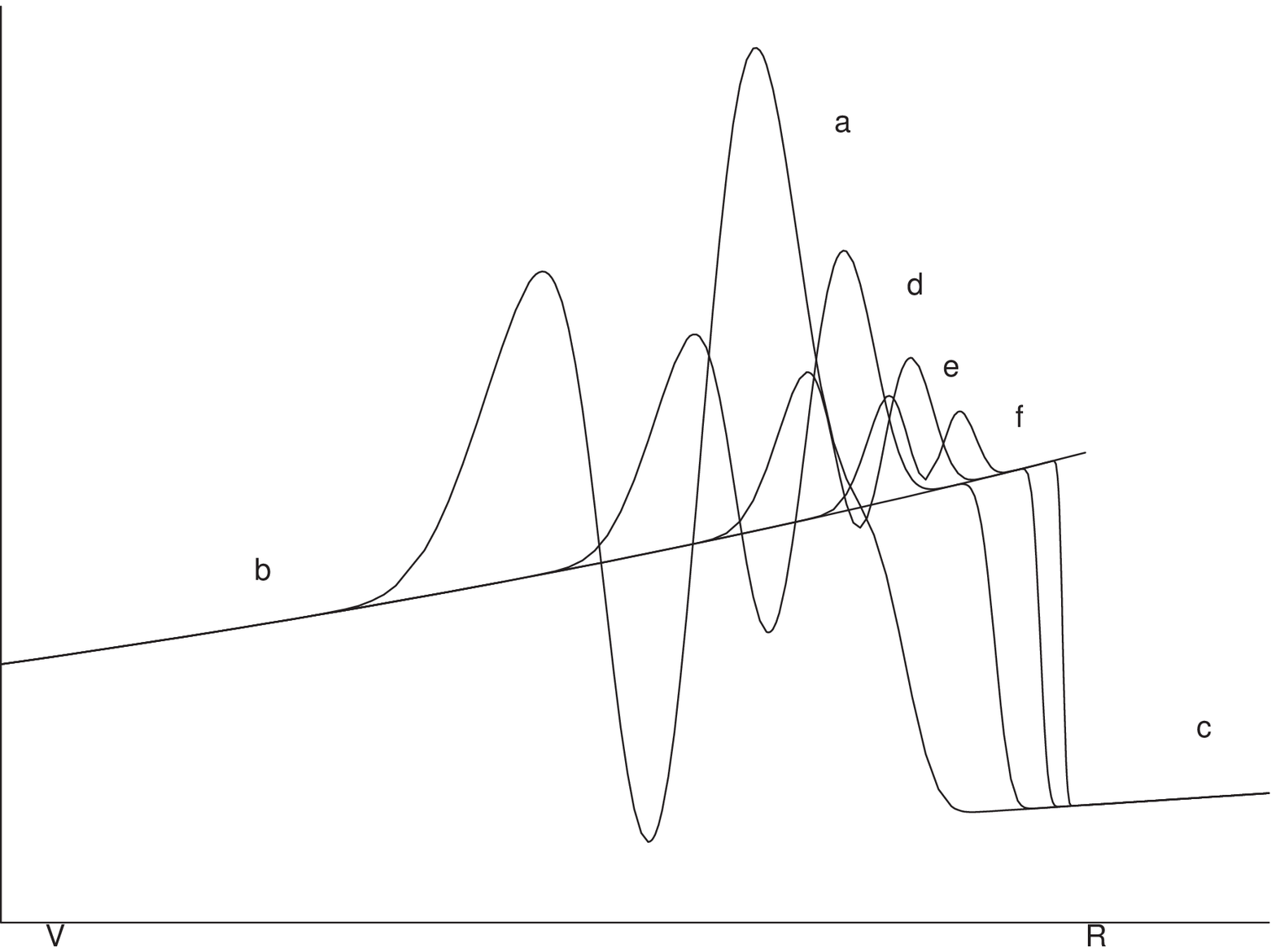,angle=0,width=.45\textwidth}}
    \caption{\label{2nd} Interpolating between $b''$ and $c''$}
    \end{huge}
    \end{center}
\end{figure}

\begin{figure}
    \begin{huge}
    \psfrag{E}{$a_\e'$}
    \psfrag{B}{$b'$}
    \psfrag{C}{$c'$}
    \psfrag{G}{$c'(R)-\eta_\e(R)$}
    \psfrag{R}{\raisebox{-2ex}{$R$}}
    \psfrag{a}{$a_\e$}
    \psfrag{b}{$b$}
    \psfrag{c}{$c$}
    \psfrag{k}{$\k_\e$}
    \psfrag{r}{$r$}
    \psfrag{D}{$\eta_\e$}
    \psfrag{U}{$R-\e$}
    \psfrag{e}{$|c(R)-\k_\e(R)|$}
    \psfrag{T}{\raisebox{-2ex}{$R-\d(\e)$}}
    \psfrag{F}{$\k_\e'$}
    \subfigure{\epsfig{figure=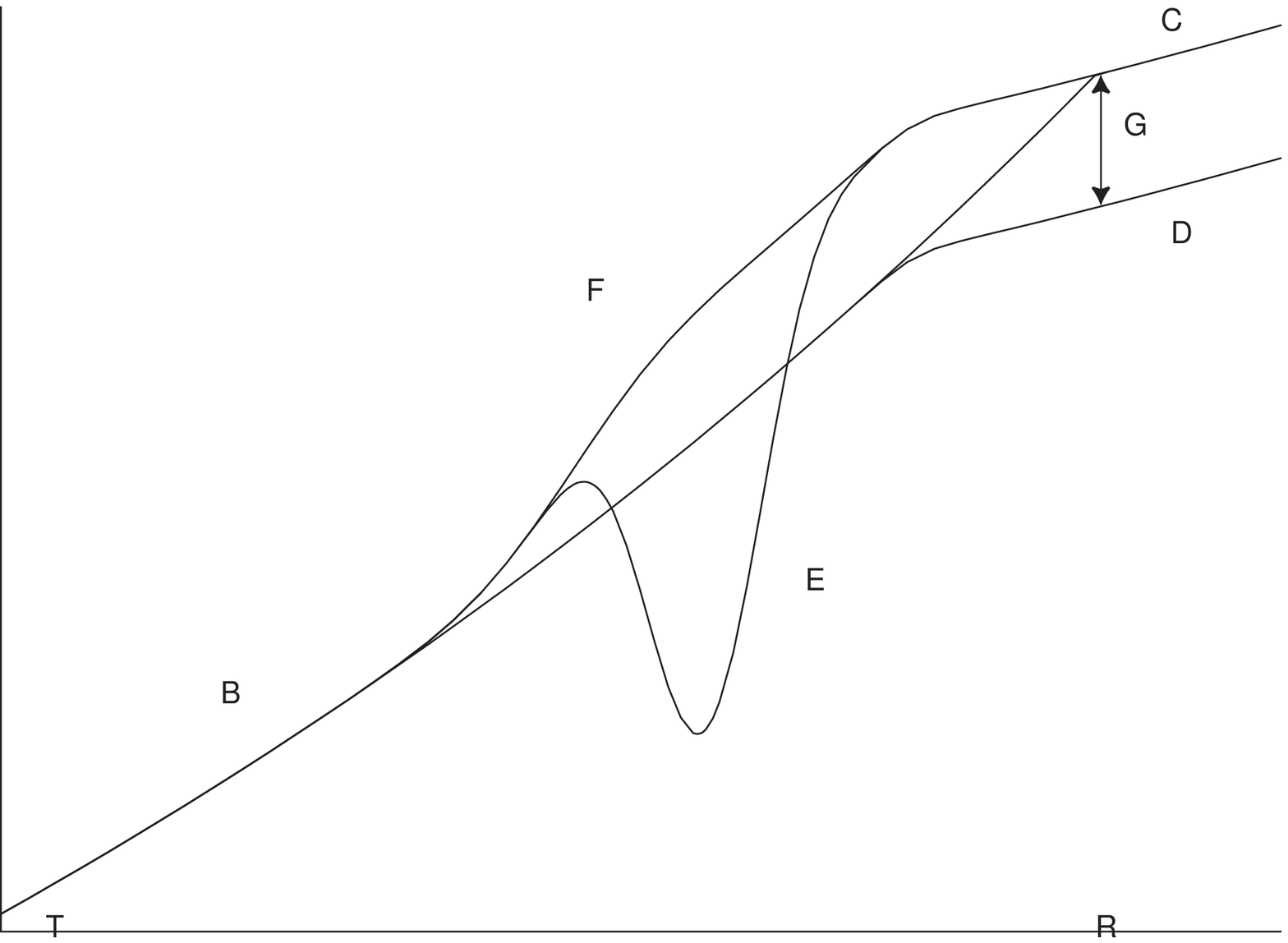,angle=0,width=.45\textwidth}}\quad
    \subfigure{\epsfig{figure=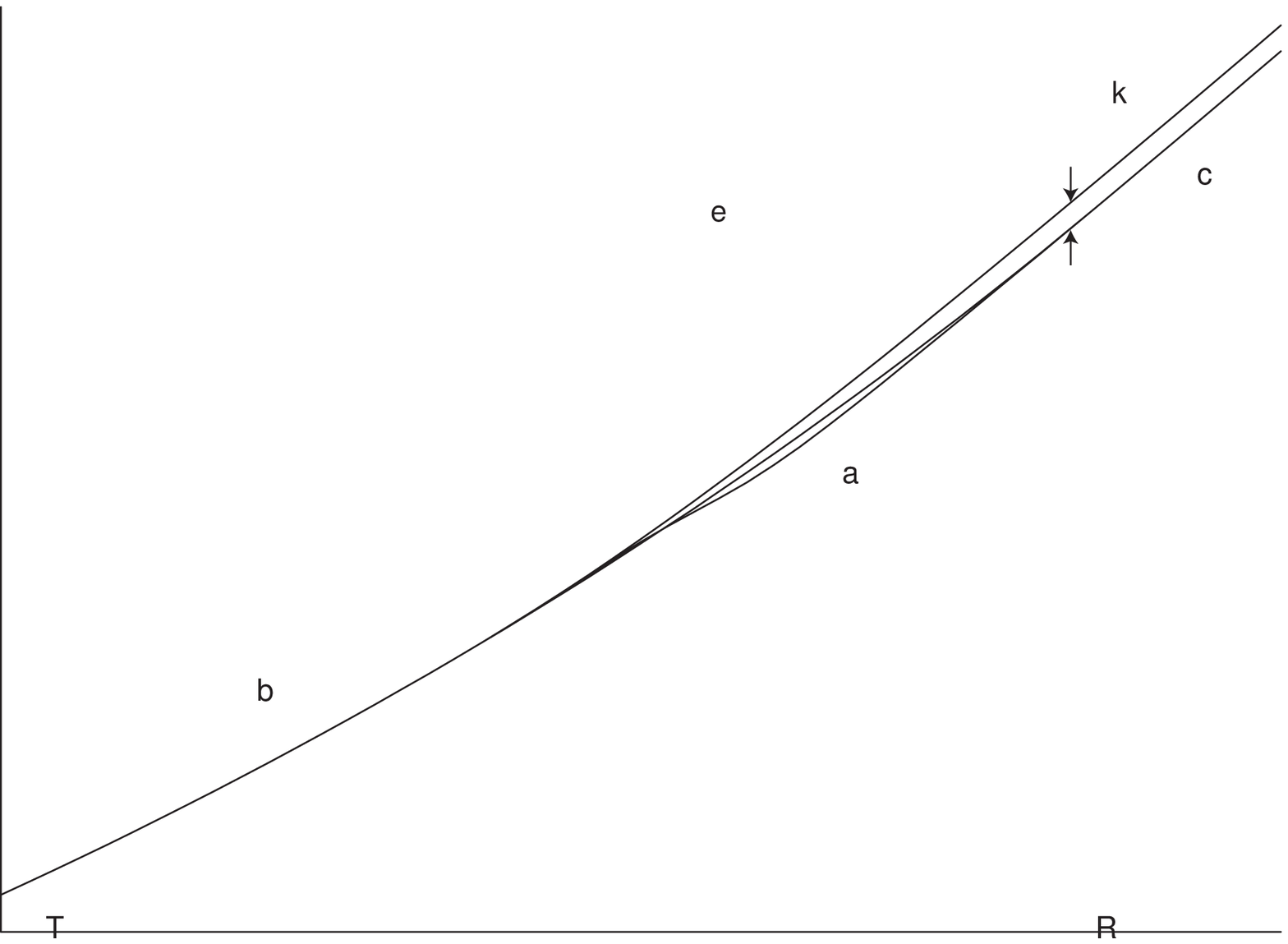,angle=0,width=.45\textwidth}}
    \caption{\label{deriv} Smoothing $a'$ and $a$}
    \end{huge}
\end{figure}

Continuing with the proof of lemma \ref{sequence}, one may see by
integrating that under the conditions of lemma \ref{smooth}, we have
$a_\e(r)\to a(r)$ and $a_\e'(r)\to a'(r)$ uniformly as $\e\to 0^+$
(figure \ref{deriv}). For the functions $f(r)$ and $g(r)$ defining
$\rho$, we can use lemma \ref{smooth} to construct a sequence of
functions $f_\e(r)$ and $g_\e(r)$ such that $f_\e'(r)\to f'(r)$,
$f_\e(r)\to f(r)$ as $\e\to 0^+$, and $\underset{\e\to
0^+}{\lim}\underset{R-\d(\e)<r<R}{\sup} f_\e''(r) = f''(R)$. Similarly
$g_\e'(r)\to g'(r)$, $g_\e(r)\to g(r)$ and $\underset{\e\to
0^+}{\lim}\underset{R-\d(\e)<r<R}{\sup} g_\e''(r) = \cosh(R)$, where
$\d(\e)$ comes from the construction in lemma \ref{smooth}. Moreover,
for $\e$ small enough, $g''(r)$ and $f''(r)>0$ for all $r$, so that
$\rho_\e$ is negatively curved for small $\e$ (figure
\ref{2nd}(b)). One may use these properties to show that $\sup
f_\e''(r)/f_\e(r) \to \coth^2 R$, $\sup g_\e''(r)/g_\e(r) \to 1$, and
$\sup f_\e'(r)g_\e'(r)/(f_\e(r) g_\e(r)) \to 1$ as $\e\to 0^+$. From
this, we deduce that $\lim_{\e\to 0^+} k_\e = (1+\coth^2 R)/2 =\coth R
\coth 2R$.
\end{proof}

\section{Applications}
First, we state a theorem of Gabai, Meyerhoff, and Thurston.

\begin{theorem}{\rm\cite[4.11]{GMT}}\label{gmt}\qua
Let $\g$ be a minimal length geodesic in a closed orientable
hyperbolic 3--manifold $M$, of length $l $ and tube radius  $R$.
Either
\begin{enumerate}
\item
$R>(\ln 3)/2$, or
\item
$1.0953/2 > R > 1.0591/2$ and $l> 1.059$, or
\item
$R=.8314.../2$ and $M=Vol3$, the closed manifold of third smallest
volume in the Snappea census \cite{W}, which has volume $=
1.0149...$.
\end{enumerate}
\end{theorem}

The following corollary is of interest in controlling the geometry
of a minimal volume hyperbolic 3--manifold. This improves on the
estimate of theorem 1.1 of Gabai, Meyerhoff, and Milley in
\cite{GMM}, and theorem 4.7 of Przeworski \cite{Pr}. This estimate
may be useful in a computational approach to find the minimal
volume orientable hyperbolic 3--manifold proposed by Gabai,
Meyerhoff, and Milley \cite{GMM} by extending the computations and
method of \cite{GMT}.

\begin{corollary} \label{minvol}
The minimal volume orientable hyperbolic 3--manifold $M$ has $\vol
(M)> .32$, and the minimal length geodesic  $\g$ in $M$ has tube
radius $R<.956$.
\end{corollary}
\proof
We may assume that $M$ is closed, since if $M$ had a cusp, then by
Cao and Meyerhoff's result that a cusped oriented manifold has
volume $>2.0298$ \cite{CM}, which is greater than the volume of
the Weeks manifold, $M$ would not have minimal volume. We may
assume that $R>(\ln 3)/2$, since in cases 2 and 3 of theorem
\ref{gmt}, the volume is $>1.01$, which is greater than the volume
of the Weeks manifold, as observed in \cite{GMT}.

 By a result
of Cao and Meyerhoff \cite{CM}, $\vol(M_\g)>2.0298$. Then by theorem \ref{bound}, we have
$\vol(M)>\vol(M_\g) / (\coth^{\frac52}(\ln 3/2)\coth^{\frac12}(\ln 3))>.32$.

We also have $$2.0298 < \vol(M_\g) \leq (\coth R)^{\frac52}(\coth
2R)^{\frac12} \vol(M) < (\coth R)^{\frac52}(\coth
2R)^{\frac12}.943.$$ Thus, we get the upper bound $R<R_0 < .956$,
where $$2.0298 =(\coth R_0)^{\frac52}(\coth 2R_0)^{\frac12}.943.\eqno{\qed}$$

\begin{figure}[htb]
    \begin{center}
    \begin{LARGE}
    \psfrag{D}{}
    \psfrag{p3}{\centering{geodesics ordered by length}}
    \psfrag{1}{$1$}
    \psfrag{10}{10}
    \psfrag{20}{20}
    \psfrag{30}{30}
    \psfrag{40}{40}
    \psfrag{8}{8}
    \psfrag{12}{12}
    \psfrag{14}{14}
    \psfrag{16}{16}
    \psfrag{2}{2}
    \psfrag{3}{3}
    \psfrag{4}{4}
    \psfrag{0}{0}
    \psfrag{6}{6}
    \psfrag{p2}{$\circ$ $\vol(M_\g)$}
    \psfrag{p1}{ $\square$ $\vol(M)+\pi l(\g) $}
    \epsfig{figure=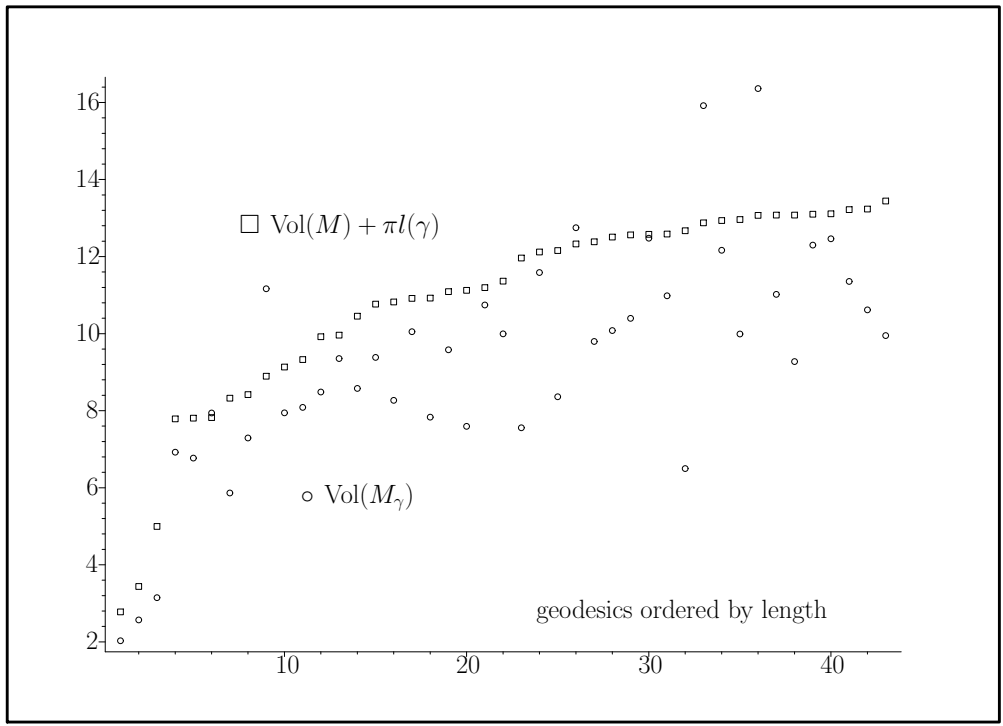,angle=-90,width=\textwidth}
    \caption{\label{drilling2} Bridgeman's conjecture for drilling geodesics in Weeks' manifold}
    \end{LARGE}
    \end{center}
\end{figure}
\begin{figure}[htb]
    \begin{center}
    \begin{LARGE}
    \psfrag{D}{}
    \psfrag{p3}{\centering{geodesics ordered by length}}
    \psfrag{1}{$1$}
    \psfrag{2}{2}
    \psfrag{3}{3}
    \psfrag{4}{4}
    \psfrag{5}{5}
    \psfrag{0}{0}
    \psfrag{6}{6}
    \psfrag{10}{10}
    \psfrag{20}{20}
    \psfrag{30}{30}
    \psfrag{40}{40}
    \psfrag{p2}{$\diamond$ $\log_{10}\vol(M_\g)$}
    \psfrag{p1}{ $\square$ $\log_{10}\coth^{5/2}R\coth^{1/2}2R\ \vol(M)$}
    \epsfig{figure=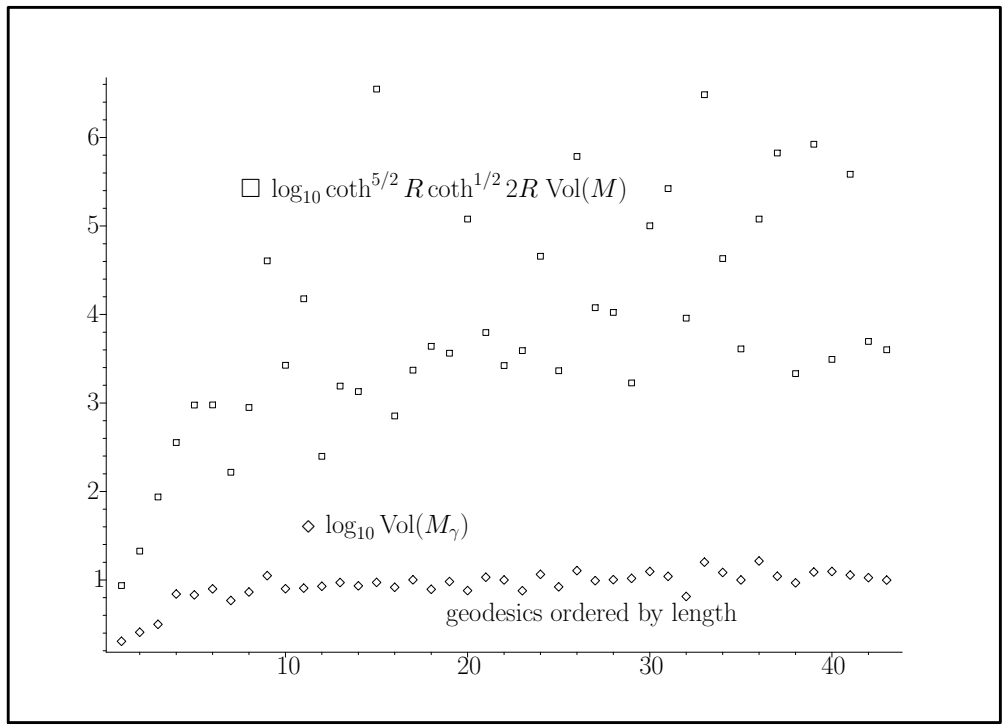,angle=-90,width=\textwidth}
    \caption{\label{drilling} Log plot of the estimate for drilling geodesics in Weeks' manifold}
    \end{LARGE}
    \end{center}
\end{figure}

{\bf Remark}\qua Przeworski has used theorem \ref{bound} and a result
of Marshall and Martin \cite{MM} to show that the minimal length
geodesic in a 3--manifold of minimal volume is $\geq .184$ and has
tube radius $R<.946$ (Proposition 7.4 \cite{Pr2}). Przeworski points out
that this estimate could be improved slightly. He has also used
theorem \ref{bound} and some packing arguments to improve the
lower bound in theorem \ref{minvol} to $.3315$ (Proposition 5.4,
\cite{Pr3}).

It would be interesting to understand exactly how much the volume
can increase when one drills out a geodesic. Bridgeman has
conjectured \cite{Bri} that if one drills an embedded geodesic
$\g$ of length $l$ out of $M$, then $\vol(M_\g)\leq \vol(M)+\pi
l$. Bridgeman proved this for many special examples. We have used
Oliver Goodman's program {\it tube} \cite{Go1} to give
computational evidence that this conjecture is false in general.
Figure \ref{drilling2} gives data obtained from {\it tube}, by
drilling many of the shortest geodesics from the Weeks manifold.
There appear to be five cases where Bridgeman's conjecture fails.
The estimate given in theorem \ref{bound} is not at all close to
being sharp, especially when $R$ is small, see figure 
\ref{drilling}.

\end{document}